\documentclass[12pt,a4paper]{article}
\usepackage{amssymb,amsfonts,amsmath}
\usepackage{theorem}

\newcommand{\rr}{\mathbb R}
\newcommand{\zz}{\mathbb Z}
\newcommand{\nn}{\mathbb N}

\theoremstyle{plain}
\theorembodyfont{\sl}

\newtheorem{theorem}{Theorem}

\newtheorem{prop}{Proposition}

\newtheorem{lemma}{Lemma}

\theorembodyfont{\rm}

\newtheorem{defn}{Definition}

\newcommand{\Cl}{\mathop{\rm Cl}\nolimits}
\newcommand{\Int}{\mathop{\rm Int}\nolimits}
\newcommand{\diam}{\mathop{\rm diam}\nolimits}
\newcommand{\dist}%
{\mathchoice{\mbox{\rm d}}{\mbox{\rm d}}%
{\mbox{\scriptsize\rm d}}{\mbox{\tiny\rm d}}}

\newenvironment{proof}[1][\hspace{-0.8ex}]%
{\noindent%
{\bf Proof #1.} {} }%
{$\blacksquare$\medskip}

\newcommand{\from}[1]{{\rm \cite{#1}}{ }}

\begin{document}

\title{On the imbedding of a finite family of closed disks into
$\rr^{2}$ or $S^{2}$}
\author{Eugene Polulyakh}
\date{\makeatletter
        Institute of mathematics,\\
        National Acad. of Sci., Ukraine\\
        \medskip{\sl e-mail polulyah@imath.kiev.ua}
        \makeatother}

\maketitle

\begin{abstract}
Let $\{V_{i}\}_{i=1}^{n}$ be a finite family of closed subsets of a
plane $\rr^{2}$ or a sphere $S^{2}$, each homeomorphic to the
two-dimensional disk. In this paper we discuss the question how
the boundary of connected components of a complement $\rr^{2}
\setminus \bigcup_{i=1}^{n} V_{i}$ (accordingly, $S^{2} \setminus
\bigcup_{i=1}^{n} V_{i}$) is arranged.

It appears, if a set $\bigcup_{i=1}^{n} \Int V_{i}$ is connected, that
the boundary $\partial W$ of every connected component $W$ of the
set $\rr^{2} \setminus \bigcup_{i=1}^{n} V_{i}$ (accordingly,
$S^{2} \setminus \bigcup_{i=1}^{n} V_{i}$) is homeomorphic to a
circle.
\end{abstract}

Let $U \in \rr^{2}$ be an open area (subset of a plane,
homeomorphic to the two-dimensional disk). One of the classical
problems of complex analysis is the question of a possibility of an
extention of conformal mapping defined in $U$ out of this area. The
answer to this question is tightly connected with the structure of
the boundary $\partial U$ of $U$ and depends on how much the
closure $\Cl U$ differs from the closed two-dimensional disk. As a
rule, it is known only the local information about a structure of
the set $\partial U$ (accessibility of points of the boundary from
area $U$ and so on).

In works~\cite{Disk,preprint} the criterion is given for a compact
subset of a plane to be homeomorphic to the closed two-dimensional
disk, which uses only local information about the boundary of this
set (see theorem~\ref{theoremD} below). This criterion enables to
investigate the problems connected to a mutual disposition of closed
disks on a plane.

Let $\{V_{i}\}_{i=1}^{n}$ be a finite family of closed subsets of a
plane $\rr^{2}$ or a sphere $S^{2}$, each homeomorphic to the
two-dimensional disk. In this paper we discuss the question how
the boundary of connected components of a complement $\rr^{2}
\setminus \bigcup_{i=1}^{n} V_{i}$ (accordingly, $S^{2} \setminus
\bigcup_{i=1}^{n} V_{i}$) is arranged.

It appears, if a set $\bigcup_{i=1}^{n} \Int V_{i}$ is connected, that
the boundary $\partial W$ of every connected component $W$ of the
set $\rr^{2} \setminus \bigcup_{i=1}^{n} V_{i}$ (accordingly,
$S^{2} \setminus \bigcup_{i=1}^{n} V_{i}$) is homeomorphic to a
circle (see. theorems~\ref{theorem_1}, \ref{theorem_2} below).

\begin{theorem}\label{theorem_1}
Let $V_{1}, \ldots, V_{n}$ be a finite collection of the closed
subsets of $\rr^{2}$, each homeomorphic to the two-dimensional
disk. Suppose the set $\bigcup_{i=1}^{n} \Int V_{i}$ is connected.

Let $W$ be the unlimited connected component of the set $\rr^{2}
\setminus \bigcup_{i=1}^{n} V_{i}$.

Then the set $\rr^{2} \setminus W$ is homeomorphic to the closed
two-dimensional disk.
\end{theorem}

\begin{theorem}\label{theorem_2}
Let $V_{1}, \ldots, V_{n}$ be a finite collection of the closed
subsets of $S^{2}$, each homeomorphic to the two-dimensional
disk. Suppose the set $\bigcup_{i=1}^{n} \Int V_{i}$ is connected
and $S^{2} \setminus \bigcup_{i=1}^{n} V_{i} \neq \emptyset$.

Let $W$ be a connected component of the set $S^{2} \setminus
\bigcup_{i=1}^{n} V_{i}$. Then the set $\Cl W$ is homeomorphic to
the closed two-dimensional disk.
\end{theorem}

The following definitions and statements will be useful for us in what
follows.

\begin{defn}\from{Tzishang}
Let $D$ be an open set. The point $x \in \partial D$ is called
\emph{accessible} from $D$ if there exists a continuous injective
mapping $\varphi : I \rightarrow \Cl D$, such that $\varphi (1) = x$ and
$\varphi ([0,1)) \subset \Int D$ (this map is named \emph{a cut}).
\end{defn}

\begin{defn}\from{Tzishang}
Let $E$ be a subset of a topological space $X$ and $a \in X$ be a point.
The set $E$ is called \emph{locally arcwise connected} in $a$,
if any neighbourhood $U$ of $a$ contains such neighbourhood $V$ of
$a$ that any two points from $V \cap E$ can be connected by a path
in $U \cap E$.
\label{defnD}
\end{defn}

\begin{prop}\from{Tzishang}
Let $D$ be an area with a nonempty interior in $\rr^{2}$ or $S^{2}$.
If $D$ is locally arcwise connected in a point $a \in \partial D$ then
$a$ is accessible from $D$.
\label{propD}
\end{prop}

\begin{theorem}\from{Disk,preprint}
Let $D$ be a compact subset of a plane $\rr^{2}$ with a nonempty
interior. Then $D$ is homeomorphic to the closed two-dimensional disk
if and only if the following conditions holds:
\begin {itemize}
        \item [1)] the set $\Int D$ is connected;
        \item [2)] the set $\rr^{2} \setminus D$ is connected;
        \item [3)] any point $x \in \partial D$ is accessible from $\Int D$;
        \item [4)] any point $x \in \partial D$ is accessible from $\rr^{2}
                  \setminus D$.
\end{itemize}
\label{theoremD}
\end{theorem}

\begin{theorem}[Sh\"{o}nflies]\label{theoremS}\from{Tzishang}
Let $\gamma$ be a simple closed curve in $S^{2}$ (respectively,
in $\rr^{2}$). There exists a homeomorphism $f$ of $S^{2}$ onto
itself (respectively, of $\rr^{2}$ onto itself) mapping the curve
$\gamma$ onto the unit circle.
\end{theorem}

\bigskip
\begin{proof}[of theorem~\ref{theorem_1}]
Let us show, that the compact set $D = \rr^{2} \setminus W$ complies
with the conditions of theorem~\ref{theoremD}. We will divide our
argument into several steps.
\medskip

{\bf 1. }
Since $\partial D \subset \bigcup_{i=1}^{n} \partial V_{i}$ then
for any $x \in \partial D$ we can find $i \in \{1, \ldots, n \}$
such that $x \in \partial V_{i}$. Theorem~\ref{theoremD} states
that the point $x$ is accessible from $\Int V_{i}$. Hence $x$
is accessible from $\Int D$ because $\Int V_{i} \subset \Int D$.
\medskip

{\bf 2. }
Let us show, that any point $a \in \partial D$ is accessible from
$W = \rr^{2} \setminus D$.

Without loss of generality we can assume that the origin of
coordinates lies in $\Int D$.

We fix $a \in \partial D$. The set of all points accessible from $W$
is dense in $\partial W = \partial D$ \cite{Tzishang}, therefore
there exists a point $x_{0} \in \partial D$ accessible from $W$ which
do not coinside with $a$.

All compact subsets of $\rr^{n}$, $n \in \nn$, are known to be
limited. Therefore there exists $R > 0$ such that
$$
\bigcup_{i=1}^{n} V_{i} \subset
\left\{ x \in \rr^{2} \,|\, \dist(0,x) < R \right\} \,.
$$

We fix a point $x' \in W$ which meets an equality $|x'| = R$.
It is known (see \cite{Tzishang}) that there exists a cut
$$
\gamma_{0} : I \rightarrow \rr^{2} \,,
$$
$\gamma_{0}(0) = x_{0}$, $\gamma_{0}(1) = x'$, $\gamma_{0}((0,1])
\subset W$.

Let
$$
\tau = \min \{\, t \in I \,|\, |\gamma_{0}(t)| = R \, \} \,.
$$
According to the conditions of theorem $\tau > 0$. Let
$\gamma_{0}(\tau) = x''$. Denote a polar angle of $x''$ by $\varphi$.

Consider continuous injective mapping
$$
\gamma_{1} : \rr_{+} \rightarrow \rr^{2} \,,
$$
$$
\gamma_{1}(t) = \left\{
\begin{array}{lrr}
        \gamma_{0}(t) \quad & \mbox{when} & t \in [0, \tau) \,, \cr
        (\varphi, R+t-\tau) & \mbox{when} & t \in [\tau, + \infty) \,. \cr
\end {array}
\right.
$$
This map is an imbedding of $\rr_{+}$ into $\rr^{2}$, moreover
$\gamma_{1}(0) \in \partial W$, $\gamma_{1}(\rr_{+} \setminus \{0\})
\subset W$.
\medskip

{\bf 2.1. }
Let us show, that the open set $W \setminus \gamma_{1}(\rr_{+})$
is connected.

Consider an involution
$$
f : \rr^{2} \setminus \{0\} \rightarrow \rr^{2} \setminus \{0\}
\,,
$$
$$
f (r, \varphi) = (r^{-1}, \varphi) \,.
$$
This map is known to be a homeomorphism. Under the action of $f$ the area
$W$ will pass to an open connected set $\widetilde{W} = f(W)$. Mark that
the origin of the coordinates is an isolated point of the boundary $\partial
\widetilde{W}$ because
$$
\{\,(r, \varphi) \in \rr^{2} \,|\, r > R \, \} \subset W \,,
$$
$$
\{\,(r, \varphi) \in \rr^{2} \,|\, 0 < r < R^{-1} \, \} \subset f(W) \,.
$$
Therefore, $\widetilde{W}_{0} = \widetilde{W} \cup \{0\}$ appeares to be
the open connected set and the map
$$
\widetilde{\gamma} : I \rightarrow \rr^{2} \,,
$$
$$
\widetilde{\gamma}(t) = \left\{
\begin{array}{llr}
        f \circ \gamma_{1}(t^{-1}-1) \quad & \mbox{when} & t \in (0, 1] \,, \cr
        0 & \mbox{for} & t = 0 \,.
\end{array}
\right.
$$
is a cut of the set $\widetilde{W}_{0}$. Moreover $\widetilde{W}_{0}
\setminus \widetilde{\gamma}(I) = f(W \setminus \gamma_{1} (\rr_{+}))$.

So, for a proof of connectivity of the set $W \setminus \gamma_{1}(\rr_{+})$
it is sufficient to check the validity of the following statement.

\begin{lemma}\label{lemma_1}
Let $U \subset \rr^{2}$ be an open connected set, point $z \in \partial
U$ be accessible from $U$, $\alpha : I \rightarrow \rr^{2}$ be a cut
of $U$ with the end in $z$ (a continuous injective mapping such that
$\alpha (0) = z$ and $\alpha ((0, 1]) \subset U$).

Then the set $U \setminus \alpha (I)$ is connected.
\end{lemma}

Let us prove this statement.
Let $y = \alpha (t)$ for some $t > 0$. According to propositions
6.4.6 and 6.5.1 from \cite{Tzishang} there exists a homeomorphism
$h$ of $\rr^{2}$ onto $\rr^{2}$, such that the map
$$
h \circ \alpha = \widetilde{\alpha} : I \rightarrow \rr^{2}
$$
complies the relation
$$
\widetilde{\alpha}(t) = (t, 0) \in \rr_{+} \times \{0\} \subset
\rr^{2} \,, \quad t \in I \,.
$$

Since $\alpha([t, 1]) \subset U$, there exists an $\varepsilon > 0$
such that
$$
\widetilde{U}_{t} = \{\, x \in \rr^{2} \,|\, \dist (x,
\widetilde{\alpha} ([t, 1])) < \varepsilon \, \} \subset h(U) \,.
$$

Obviously, $U_{t} = h^{-1} (\widetilde{U}_{t})$ is a neighbourhood
of a point $\alpha(t)$ in $U$ and the set $U_{t} \setminus
\alpha(I)$ is connected. Besides, a set $(U_{t_{1}} \cap U_{t_{2}})
\setminus \alpha (I)$ is not empty for any $t_{1}$, $t_{2} \in (0,
1]$.

Therefore
$$
\bigcup_{t \in (0, 1]} U_{t}
$$
is a connected open neighbourhood of a set $\alpha (I)$ in $U$,
hence $U \setminus \alpha (I)$ is a connected set.
$\square$
\medskip

So, the set $W \setminus \gamma_{1} (\rr_{+})$ is connected.
\medskip

{\bf 2.2. }
Select a point $x_{i} \in \partial V_{i}$, $x_{i} \ne a$ for each
$i \in \{1, \ldots, n \}$. The set
$$
\bigcup_{i=1}^{n} \Int V_{i}
$$
is connected by the condition of theorem and the point $x_{i}$ is
accessible from $\Int V_{i}$ for any $i$. Therefore we can find a
continuous map
$$
\beta : [1, n+1] \rightarrow \bigcup_{i=1}^{n} V_{i}
$$
which meets the following conditions
$$
\beta ((i, i+1)) \subset \bigcup_{i=1}^{n} \Int V_{i} \subset \Int
D \,, \quad i=1, \ldots, n \,;
$$
$$
\beta (i) = x_{i} \,, \quad i=1, \ldots, n \,;
\quad \beta (n+1) = x_{0} \,.
$$

Consider a continuous map
$$
\gamma: \rr_{+} \rightarrow \rr^{2} \,,
$$
$$
\gamma (t) = \left\{
\begin{array}{llr}
        \beta (t+1) \quad & \mbox {for} & t \in [0, n) \,, \cr
        \gamma_{1} (t-n) & \mbox{for} & t \in [n, + \infty) \,. \cr
\end{array}
\right.
$$
Since the relations
$$
\gamma (\rr_{+}) \subset \left(\beta ([1, n+1]) \cup \gamma_{0} (I)
\cup \gamma_{1} ([\tau, + \infty)) \right) \,,
$$
$$
\gamma_{1} ([\tau, + \infty)) \subset \{\, z \in \rr^{2} \,|\,
\dist (z, 0) \ge R \, \} \,,
$$
$$
a \in \partial D \subset \{\, z \in \rr^{2} \,|\, \dist (z, 0) < R \, \}
$$
hold and a compact set $\beta ([1, n+1]) \cup \gamma_{0}(I)$ does not
contain a point $x$ on a construction, there exists
$\varepsilon_{0} > 0$ which complies the inequality
$$
\dist (a, z) > \varepsilon_{0} \quad \mbox{for all{ }} z \in \gamma
(\rr_{+}) \,.
$$

Now we are ready for proof of local linear connectivity of the area
$W$ in the point $a \in \partial W = \partial D$.

\medskip
{\bf 2.3. }
Let $U$ be a curtain neighbourhood of the point $a$. Find $
\varepsilon > 0$ which meets the conditions
$$
U_{\varepsilon} (a) = \{\, x \in \rr^{2} \,|\, \dist (a, x) <
\varepsilon \, \} \subset U \,, \quad U_{\varepsilon} (a) \cap
\gamma (\rr_{+}) = \emptyset \,.
$$

Fix imbeddings
$$
f_{i} : S^{1} \rightarrow \partial V_{i} \,, \quad i = 1, \ldots, n
\,.
$$
Here $S^{1} = \{\,(r, \varphi) \in \rr^{2} \,|\, r=1 \, \}$. The
metric on $S^{1}$ we shall define as follows:
$$
\dist_{s} \left((1, \varphi_{1}) ,\, (1, \varphi_{2}) \right) =
\min_{k \in \zz} | \, \varphi_{1} - \varphi_{2} + 2 \pi k \, | \,.
$$

Mark that maps $f_{i}$, $i=1, \ldots, n$ are uniformly continuous.

Fix $\delta_{1} > 0$ such that an inequality $\dist_{s} (\tau_{1},
\tau_{2}) < \delta_{1}$ implies
$$
\dist (f_{i}(\tau_{1}),\, f_{i}(\tau_{2})) < \min \left(
\varepsilon_{0} / 2 ,\, \varepsilon / 3 \right)
$$
for any $i=1, \ldots, n$ and $\tau_{1}$, $\tau_{2} \in S^{1}$.

Find also $\delta_{2} > 0$, such that $\dist (z_{1},
z_{2}) < \delta_{2}$ has as a consequence an inequality
$$
\dist_{s} (f_{i}^{-1} (z_{1}),\, f_{i}^{-1} (z_{2})) < \min \left(
\delta_{1} / 2,\, \pi / 4 \right)
$$
for every $i=1, \ldots, n $ and any $z_{1}$, $z_{2} \in \partial
V_{i}$.

Assume $\delta = \min \left(\delta_{2} / 2,\, \varepsilon / 3
\right)$.
\medskip

{\bf 2.4. }
Let us show, that for any $a_{1}$, $a_{2} \in U_{\delta} (a) \cap
W$ there exists a continuous map $g : I \rightarrow U_{\varepsilon}
(a) \cap W$ such that $g(0) = a_{1}$, $g(1) = a_{2}$.

The inequality $\dist (z_{1}, z_{2}) < \delta_{2}$ is fulfilled for
all $z_{1}$, $z_{2} \in U_{\delta} (a)$, hence
$$
\dist \left(f^{-1} (\partial V_{i} \cap U_{\delta} (a)) \right) <
\min \left( \delta_{1} / 2,\, \pi / 4 \right)
$$
for every $i \in \{1, \ldots, n \}$ and in the case $\partial V_{i}
\cap U_{\delta}(a) \neq \emptyset$ the circle $S^{1}$ could be
decomposed into two not intersecting intervals $J_{i}'$ and
$J_{i}''$ with common endpoints in such a way that the following
relations are fulfilled
$$
f_{i}^{-1} (\partial V_{i} \cap U_{\delta} (a)) \subset J_{i}' \,,
$$
$$
\diam (J_{i}') = \max_{t_{1}, t_{2} \in J_{i}'} \dist_{s} (t_{1},
t_{2}) < \min \left( \delta_{1},\, \pi / 2 \right) \,.
$$
In the case $\partial V_{i} \cap U_{\delta}(a) = \emptyset$ set
$J_{i}'' = S^{1}$, $J_{i}' = \emptyset$.

Therefore, $f_{i} (J_{i}'') \cap U_{\delta} (a) = \emptyset$ and
$f_{i}(J_{i}') \subset U_{2 \varepsilon / 3}(a)$.

\begin{lemma}\label{lemma_2}
Let $B$ be a closed disk satisfying the following conditions:
$$
\partial B \cap \left( \textstyle{\bigcup\limits_{i=1}^{n}}
\partial V_{i} \right) \subset U_{\delta} (a) \,,
$$
$$
( \partial B \setminus U_{\delta} (a)) \subset (W \setminus \gamma
(\rr_{+})) \,.
$$

Then $B \cap \partial V_{i} \subset f_{i} (J_{i}') \subset U_{2
\varepsilon / 3}(a)$, $i = 1, \ldots, n$.
\end{lemma}

On the condition of lemma $\partial B \cap f _ {i} (J _ {i}'') =
\emptyset$ for every $i = 1, \ldots, n$. Therefore, $f_{i}
(J_{i}'') \subset \Int B$ or $f_{i} (J_{i}'') \subset (\rr^{2}
\setminus B)$. By a construction $x_{i} \in f_{i} (J_{i}'')$ and
$x_{i} \in \gamma (\rr_{+}) \subset (\rr^{2} \setminus B)$, hence
$f_{i} (J_{i}'') \subset (\rr^{2} \setminus B)$, $i = 1, \ldots,
n$.
$\square$

Let $a_{1}$, $a_{2} \in (U_{\delta} (a) \cap W)$. Since $U_{\delta}
(a) \cap \gamma (\rr_{+}) = \emptyset$, then $a_{1}$, $a_{2} \in
(U_{\delta} (a) \cap (W \setminus \gamma (\rr_{+})))$. From
connectivity of the set $W \setminus \gamma (\rr_{+})$ follows,
that there exists an injective continuous map
$$
\widetilde{\mu} : I \rightarrow (W \setminus \gamma (\rr_{+})) \,,
$$
complying the equalities $\widetilde{\mu} (0) = a_{1}$, $
\widetilde{\mu} (1) = a_{2}$ (the concepts of connectivity and
linear connectivity coincide for open subsets of $\rr^{n}$).

Find smooth imbeddings
$$
\eta_{1} : S^{1} \rightarrow U_{\delta} (a) \,,
$$
$$
\eta_{2} : S^{1} \rightarrow \left(U_{\varepsilon} (a)
\setminus \Cl U_{2 \varepsilon /3} (a) \right) \,,
$$
such that the points $a_{1}$, $a_{2}$ lie inside disks bounded by
curves $\eta_{1}$, $\eta_{2}$.

It is known that an imbedding of a segment or circle into $\rr^{2}$
could be as much as desired precisely approximated by a smooth
imbedding. It is known as well that any two one-dimensional smooth
compact submanifolds of $\rr^{2}$ could be reduced to the general
position by a small perturbation fixed on their boundary.

Therefore, there exists smooth imbedding
$$
\mu : I \rightarrow W \setminus \gamma (\rr_{+}) \,, \quad
a_{1} = \mu (0) \,,\; a_{2} = \mu (1)
$$
such that the sets $\mu (I) \cap \eta_{1} (S^{1})$ and $\mu (I)
\cap \eta_{2} (S^{1})$ consist of final number of points.

For every $z \in \mu (I) \cap \eta_{2} (S^{1})$ there exist $t'$,
$t'' \in I$, $t' < t''$, which comply with the following conditions
$$
z \in \mu ((t', t'')) \,,
$$
$$
\mu (t'), \mu (t'') \in \eta_{1} (S^{1}) \,,
$$
$$
\mu ((t', t'')) \cap \eta_{1} (S^{1}) = \emptyset \,.
$$

We receive a finite family of nonintersecting intervals
$$
(t_{j, 1}, t_{j, 2}) \subset I \, \qquad j = 1, \ldots, k
$$
satisfying to relations
$$
\mu ((t_{j, 1}, t_{j, 2})) \cap \eta_{1} (S^{1}) = \emptyset \,,
\quad \mu (t_{j, 1}), \mu (t_{j, 2}) \in \eta_{1} (S^{1}) \,, \quad
j = 1, \ldots, k \,,
$$
$$
\mu \Bigl(I \bigl\backslash \textstyle{\bigcup\limits_{j=1}^{k}}
(t_{j, 1}, t_{j, 2}) \Bigr) \subset U_{\varepsilon} (a) \,.
$$

Now for each $j = 1, \ldots, k$ we fix an arc $\Theta_{j} : I
\rightarrow \eta_{1} (S^{1})$ with the endpoints $\mu (t_{j, 1})$
and $\mu (t_{j, 2})$. A set
$$
\Theta_{j} (i) \cup \mu ((t_{j, 1}, t_{j, 2}))
$$
is homeomorphic to a circle, therefore it bounds a closed disk
$B_{j} $ such that
$$
\left(\partial B_{j} \cap \textstyle{\bigcup\limits_{i=1}^{n}}
\partial V_{i} \right) \subset U_{\delta} (a) \,,
$$
$$
\left(\partial B_{j} \setminus U_{\delta} (a) \right) \subset
\left(W \setminus \gamma (\rr_{+}) \right) \,.
$$

By lemma~\ref{lemma_2} these relations has as a
consequence following inclusion
$$
\left(B_{j} \cap \textstyle{\bigcup\limits_{i=1}^{n}} \partial
V_{i} \right) \subset U_{2 \varepsilon / 3} (a) \,.
$$

Since $\eta_{2} (S^{1}) \subset (U_{\varepsilon} (a) \setminus \Cl
U_{2 \varepsilon / 3} (a))$, then
$$
B_{j} \cap \eta_{2} (S^{1}) =
\textstyle{\bigcup\limits_{s=1}^{m_{j}}} \chi_{s} \,.
$$
Here $\{\chi_{s}\}_{s=1}^{m_{j}}$ is a final family of
nonintersecting arcs of the circle $\eta_{2} (S^{1})$. In addition
$\chi_{s} \subset (W \setminus \gamma (\rr_{+}))$, $s = 1, \ldots,
m _ {j}$.

A set
$$
(\Int B_{j}) \; \bigl\backslash \; \Bigl(
\textstyle{\bigcup\limits_{s=1}^{m _ {j}}} \chi_{s} \Bigr)
$$
represents a final union of connected components homeomorphic to
the two-dimensional disk, lying either inside or outside the
closed disk limited by a circle $\eta_{2} (S^{1})$. Select that
from components, which bounds with an arc $\Theta_{j}$. Designate
by $\widetilde{B}_{j}$ a closure of this component. Obviously,
$$
\widetilde{B}_{j} \subset U_{\varepsilon} (a) \,, \quad (\partial
\widetilde{B}_{j} \setminus \Theta_{j} (I)) \subset (W \setminus
\gamma (\rr_{+})) \,.
$$

Let
$$
g_{j} : I \rightarrow (\partial \widetilde{B}_{j} \setminus
\Theta_{j} ((0, 1)))
$$
be an arc of a circle $\partial \widetilde{B}_{j}$ with the
endpoints $\mu (t_{j, 1})$, $\mu (t_{j, 2})$. As we already have
shown, it complies with the relation
$$
g_{j} (I) \subset \left(U_{\varepsilon} (a) \cap (W \setminus
\gamma (\rr_{+})) \right) \,.
$$
A continuous curve
$$
g : I \rightarrow (W \setminus \gamma (\rr_{+})) \,,
$$
$$
g(t) = \left\{
\begin{array}{lrr}
        \mu(t) & \mbox{if} & t \in \left(I \; \bigl\backslash \;
                \textstyle{\bigcup_{j=1}^{k}} (t_{j, 1}, t_{j, 2})
                \right) \,, \cr
        g_{j}((t_{j, 2} - t_{j, 1}) t + t_{j, 1}) & \mbox{if} & t
                \in (t_{j, 1}, t_{j, 2}) \,. \cr
\end{array}
\right.
$$
represents a continuous path in $U_{\varepsilon} (a) \cap (W
\setminus \gamma (\rr_{+}))$, connecting points $a_{1}$ and
$a_{2}$.

Therefore, the point $a$ is accessible from $W \setminus \gamma
(\rr_{+})$ and all the more it is accessible from $W = \rr^{2}
\setminus D$. Then each point of $\partial D$ is accessible from
$\rr^{2} \setminus D$ because of the arbitrary rule we selected
the point $a \in \partial D$.
\medskip

{\bf 3. }
The set $W$ is connected on a condition of theorem.
\medskip

{\bf 4. }
Let us show that the set $\Int D$ is connected. The set
$\bigcup_{i=1}^{n} V_{i}$ is connected since any point of
$\bigcup_{i=1}^{n} \partial V_{i}$ is accessible from a connected
set $\bigcup_{i=1}^{n} \Int V_{i}$, therefore it is sufficient to
show that the boundary $ \partial \widetilde{W}$ does not lie in
the set $\partial D$ for any connected component $\widetilde{W}$ of
the set $\rr^{2} \setminus \left(\bigcup_{i=1}^{n} V_{i} \right)$,
different from $W$.

Assume that $\partial \widetilde{W} \subset \partial D$. The set
$\partial \widetilde{W}$ divides $\rr^{2}$, consequently it has
dimension not less than one (see \cite{Gurevich}). Therefore,
we can find three different points $z_{1}$, $z_{2}$, $z_{3} \in
\partial \widetilde{W}$. Each of these points is accessible from
the connected sets $W$ and $\bigcup_{i=1}^{n} \Int V_{i}$.

There exists a continuous injective mapping (see \cite{Tzishang})
$$
\varphi : I \rightarrow \rr^{2}
$$
which satisfies the conditions
$$
\varphi(0) = z_{1} \,, \quad \varphi(1) = z_{2} \,, \quad
\varphi((0, 1)) \subset \textstyle{\bigcup\limits_{i=1}^{n}} \Int
V_{i} \,.
$$

Let $z = \varphi (1/2)$. There exists a continuous injective
mapping
$$
\widetilde{\varphi} : I \rightarrow \rr^{2} \,,
$$
$$
\widetilde{\varphi} (0) = z_{3} \,, \quad \widetilde{\varphi} (1) =
z \,, \quad \widetilde{\varphi} ((0, 1]) \subset
\textstyle{\bigcup\limits_{i=1}^{n}} \Int V_{i} \,.
$$

Let $t_{1} = \min \{\, t \in I \,|\, \widetilde{\varphi} (t) \in
\varphi (I) \, \}$. We have $t_{1} > 0$ since $z_{3} =
\widetilde{\varphi} (0) \notin \varphi (I)$. Denote $z' =
\widetilde{\varphi} (t_{1})$. Then $t_{2} \in (0, 1)$ is uniquely
defined, such that $z' = \varphi (t_{2})$.

Consider continuous injective mappings
$$
\varphi_{1} : I \rightarrow \rr^{2}, \quad \varphi_{1} (t) =
\varphi (t_{2} (1-t)) \,;
$$
$$
\varphi_{2} : I \rightarrow \rr^{2}, \quad \varphi_{2} (t) =
\varphi ((1-t_{2}) t + t_{2}) \,;
$$
$$
\varphi_{3} : I \rightarrow \rr^{2}, \quad \varphi_{3} (t) =
\widetilde{\varphi} (t_{1} (1-t)) \,.
$$
which comply with the relations
$$
\varphi_{s} (0) = z_{s}, \quad \varphi_{s} (1) = z', \quad
\varphi_{s} ((0, 1]) \subset \textstyle{\bigcup\limits_{i=1}^{n}}
\Int V_{i} \,, \quad s=1, 2, 3 \,;
$$
$$
\varphi_{s_{1}} ([0, 1)) \cap \varphi_{s_{2}} ([0, 1)) = \emptyset
\qquad \mbox{when{ }} s_{1} \ne s_{2} \,.
$$

Similarly, there exists a point $z'' \in W$ and continuous
injective mappings $\psi_{s} : I \rightarrow \rr^{2}$, $s=1,2,3$,
such that
$$
\psi_{s} (0) = z_{s}, \quad \psi_{s} (1) = z'', \quad \psi_{s} ((0,
1]) \subset W \,, \quad s=1, 2, 3 \,;
$$
$$
\psi_{s_{1}} ([0, 1)) \cap \psi_{s_{2}} ([0, 1)) = \emptyset \qquad
\mbox{when{ }} s_{1} \ne s_{2} \,.
$$

Since
$$
W \cap \Bigl(\textstyle{\bigcup\limits_{i=1}^{n}} \Int V_{i} \Bigr)
= \emptyset \,,
$$
the equality
$$
\left( \Bigl( \textstyle{\bigcup\limits_{s=1}^{3}} \varphi_{s} (I)
\Bigr) \cap \Bigl( \textstyle{\bigcup\limits_{s=1}^{3}} \psi_{s}
(I) \Bigr) \right) = \textstyle{\bigcup\limits_{s=1}^{3}} \{z_{s}
\}
$$
is valid. Therefore, everyone from the sets
$$
\varphi_{s_{1}} (I) \cup \varphi_{s_{2}} (I) \cup \psi_{s_{1}} (I)
\cup \psi_{s_{2}} (I) \,, \quad s_{1} \ne s_{2}
$$
is homeomorphic to a circle.

The set
$$
\rr^{2} \; \bigl\backslash \; \Bigl(
\textstyle{\bigcup\limits_{s=1}^{3}} (\varphi_{s} (I) \cup \psi_{s}
(I)) \Bigr)
$$
falls into the three connected components $U_{1}$, $U_{2}$,
$U_{3}$, two of which are homeomorphic to the open two-dimensional
disk and third is not limited.

As
$$
\Bigl( \textstyle{\bigcup\limits_{s=1}^{3}} (\varphi_{s} (I) \cup
\psi_{s} (I)) \Bigr) \cap \widetilde{W} = \emptyset \,,
$$
then there exists $j \in \{1, 2, 3 \}$ such that $\widetilde {W}
\subset U_{j}$. But it is impossible because everyone from the sets
$\Cl U_{s}$, $s=1, 2, 3$ contains exactly two from the points
$z_{1}$, $z_{2}$, $z_{3}$.

So, we have proved that the set $\partial \widetilde{W} \cap
\partial D$ consists not more than from two points. Therefore, $z
\in \partial \widetilde {W}$ and $\varepsilon > 0$ could be found
to comply the inclusion $U_{\varepsilon} (z) \subset \Int D$.

The set
$$
\widetilde {W} \cup U_{\varepsilon} (z) \cup \Bigl(
{\textstyle\bigcup\limits_{i=1}^{n}} \Int V_{i} \Bigr) \subset
\Int D
$$
is connected since $\partial \widetilde{W} \subset
\bigcup_{i=1}^{n} \partial V_{i}$ and the sets $\widetilde {W}$,
$U_{\varepsilon} (z)$, $\bigcup_{i=1}^{n} \Int V_{i}$ are
connected.

By virtue of arbitrariness in a choice of $\widetilde{W}$, the set
$\Int D$ is connected.

Applying to $D$ theorem~\ref{theoremD} we conclude that this set
is homeomorphic to the closed two-dimensional disk.
\end {proof}

\bigskip
\begin{proof}[of theorem~\ref{theorem_2}]
Let $in_{i} : I^{2} \rightarrow S^{2}$, $i = 1, \ldots, n$ be
the inclusion maps, $in_{i}(I^{2}) = V_{i}$.

Without loss of a generality, it is possible to assume that
a North Pole $s_{0}$ of $S^{2}$ lies in $W$.

Consider a stereographic projection
$$
f : S^{2} \setminus \{s_{0}\} \rightarrow \rr^{2} \,.
$$
As is known, this map is a homeomorphism. Since
$V_{i} \subset S^{2} \setminus \{s_{0}\}$, $i = 1, \ldots, n$ and
the set $S^{2} \setminus \{s_{0}\} $ is open in $S^{2}$,
the compositions
$$
In_{i} = f \circ in_{i} : I^{2} \rightarrow \rr^{2} \,,
\quad i = 1, \ldots, n
$$
are continuous and are one-to-one. The set $I^{2}$ is compact,
therefore maps $In_{i}$, $i = 1, \ldots, n$ are
imbeddings. Sign $\widehat{V}_{i} = f(V_{i}) =
In_{i}(I^{2})$, $i = 1, \ldots, n$.

From a mutual uniqueness of map $f$ follows that
$$
f(\bigcup_{i=1}^{n} \Int V_{i}) = \bigcup_{i=1}^{n} f (\Int V_{i})
= \bigcup_{i=1}^{n} \Int \widehat{V}_{i} \,.
$$
The set $\bigcup_{i=1}^{n} \Int \widehat{V}_{i}$
is connected as an image of a connected set at a continuous map.

So, family $\widehat{V}_{1}, \ldots, \widehat{V}_{n}$
satisfies to conditions of theorem~\ref{theorem_1}.

Consider an open set $W' = W \setminus \{s_{0}\} \subset
S^{2}$. It is easy to see that $\partial W' = \partial W \cup
\{s_{0}\}$ and $s_{0}$ is an isolated point of the boundary
of $W'$.

Denote $\widehat{W} = f(W') \subset \rr^{2}$. Obviously, $\widehat{W}$
is the unique unlimited connected component of a set $\rr^{2}
\setminus \bigcup_{i=1}^{n} \widehat{V}_{i}$. Applying
theorem~\ref{theorem_1}, we conclude that a set $\rr^{2} \setminus
\widehat{W}$ is homeomorphic to the closed two-dimensional disk, and it's
boundary $\partial (\rr^{2} \setminus \widehat{W}) = \partial \widehat{W}$
is homeomorphic to a circle $S^{1}$. From this immediately follows, that
the set $\partial W = f^{-1} (\partial \widehat{W})$ of the limit points
of $W$ is homeomorphic to a circle.

From theorem~\ref{theoremS} it immediately follows that the set
$\partial W$ divides $S^{2}$ into two opened connected components
and for each of these components it's closure is homeomorphic to
the closed two-dimensional disk. Consequently, the set $\Cl W$ is
homeomorphic to the closed two-dimensional disk.
\end{proof}

\end{document}